\documentclass[12pt]{article}
\usepackage[utf8]{inputenc}
\usepackage{amsmath}
\usepackage{lmodern,amssymb,bbm}
\usepackage{latexsym}
\usepackage{graphicx}
\usepackage{hyperref}

\usepackage[left=3cm,right=3cm,top=3cm,bottom=3cm]{geometry}
\usepackage{color}
\usepackage{cancel}
\usepackage{comment}

\definecolor{brown}{RGB}{150,100,0}
 
\definecolor{pink}{RGB}{255,0, 255}
 
\definecolor{grey}{RGB}{128,128,128}

\newtheorem{theorem}{Theorem}
\newtheorem{lemma}{Lemma}
\newtheorem{definition}{Definition}
\newtheorem{proposition}{Proposition}

\usepackage{soul}
\usepackage{comment}

\definecolor{purple}{RGB}{128,0,128}

\title{Epidemic Dynamics via Wavelet Theory 
and Machine Learning, with Applications to Covid-19\footnote{This work is  supported by Torus Actions and Belle Artificial Intelligence Corporation}}
\author{T\^o Tat Dat\footnote{corresponding author: tat-dat.to@imj-prg.fr} , Protin Frédéric,  Nguyen T.T. Hang,\\
Martel Jules\footnote{Invited fellow at Max Planck Institute, Bonn},   Nguyen Duc Thang,  Charles Piffault, Rodr\'iguez Willy, \\  Figueroa Susely,
 H\^ong V\^an  L\^e\footnote{partially supported  by GA\v CR-project 18-01953J and	 RVO: 67985840}, Wilderich Tuschmann, Nguyen Tien Zung}
\date{}

\begin{document}
\maketitle

\begin{abstract}
\noindent
We   introduce the  concept of epidemic-fitted  wavelets which comprise, in particular, as special cases
the  number $I(t)$  of infectious individuals  at time $t$ in classical SIR models and their derivatives. 
We present a novel
 method for   modelling     epidemic  dynamics 
 by a model selection  method using wavelet theory and,
 for its applications, machine learning based curve fitting techniques. 
 Our universal models     are    functions that are finite linear combinations of epidemic-fitted wavelets. 
We  apply  our  method
by modelling and forecasting, based on the John Hopkins University dataset, the spread of the current Covid-19 (SARS-CoV-2) epidemic in France, Germany, Italy and  the Czech Republic, as well as
in the US federal states New York and Florida.
\end{abstract}

\section{Introduction}


\

\noindent
The present work proposes  a  novel    method
for modelling   epidemic  dynamics by  combining
wavelet theory  and data-driven model section 
techniques in machine  learning. 

In understanding epidemic diffusion and growth rate of an infectious disease
at population level, 
the actual number of reported cases of infections 
always plays a (if not even \emph{the}) crucial role - and,
way beyond that, at least in the case of diseases afflicting human societies, directly influences government and health care system decisions and measures regarding, e.g., protection, 
containment and hospital capacities.

However, due to both the manifold practical as well as conceptual issues
involved,
a rigorous and accurate detection of this number
turns out to be a rather difficult and complex problem.

 To illustrate at least some of the theoretical difficulties 
 involved here 
 by a prominent and important case which calls the entire world now to action,  let us note that most current mathematical modelling and forecasting techniques for the spread of the Covid-19 disease are based
 on classical SIR (Susceptible - Infectious - Recovered/Removed) and SEIR 
 (Susceptible - Exposed - Infectious - Recovered/Removed) compartmental 
 epidemiological models \cite{KM,BDW2008,Wang}.
 
 Yet, with regard to predicting
 the number  of   infectious  cases $I(t)$  at time $t$,
they suffer from severe and model-inherent principal limitations:
 
All these models, as well as all their derivatives, 
are  not  suitable  to build a model for the function $I(t)$  
which is compatible with any given population. 
This is because these models 
are based on the assumption that the population is homogeneously
composed and distributed (i.e., the chance
that an arbitrary infected person will infect an
arbitrary susceptible person is taken to be constant throughout the epidemic, and, moreover, 
it is assumed that at any given time every infected person has one 
and the same constant chance  to recover). 

In real life, however, 
there are  actually many and rather diverse waves of outbreaks, 
stemming from different times or locations. 
One faces here not only drastically varying
growth rates, but also hot spots versus no-cluster
locations, infection rates depending on age or other parameters, etc.
which altogether entails that 
the homogeneity assumption approach
taken in SIR models and their variations
is over-simplified and cannot give realistic forecasts.

To overcome the draw-backs caused by homogeneity assumptions, the new approach presented in this work is based on the following idea:
we  shall decompose the growth curve of infection
numbers into  several   basic  'waves', where each  basic wave is considered as a representation 
of  the epidemic, and localised both in time and position.

This point of view naturally calls for the use of  wavelet theory.
Wavelets as such are special families of functions which
came up in the 1980s by combining older concepts from
mathematics, computer science, electrical engineering and physics,
having since then found fruitful applications in many other disciplines.
In particular, some precursor, wave-based approaches to modelling epidemic growth appeared already a long time before wavelets emerged in both deterministic and stochastic  models, 
compare, among others,  \cite{Sop,Bar1956, Bar1957, KR},
and only very recently, Krantz et al.
(compare \cite{KPR})
Moreover, the latter work has also proposed
building epidemic growth models by combining wavelet with discrete graph theory (see also below).  

In this article, we propose an approach to  epidemic dynamics
by modelling the number of daily reported cases using specially designed wavelets, called  epidemic-fitted (EF) wavelets. 

For instance, the number  $I(t)$  of infectious  individuals at time $t$  in the 
classical SIR and SEIR  models   is  an EF wavelet, see Subsection \ref{ssect:choose}. Another example of an EF wavelet is  the log-normal
one, which we will use in our Covid-19 spread forecasting applications,
see   subsections \ref{ssect:log_normal}, \ref{ssect:choose}  for   more details.

In our approach, the number  of daily reported  cases   is   the value of a  function that is a   positive linear  combination of  $N$ EF wavelets  at  the  given day.   We fix the  number 
$N$ of   summands  of   EF wavelets entering  in   our modelling  function (and in our applications $N$ is usually taken to be $3$ or $5$). 
The wavelet series coefficients themselves 
are  then obtained by machine learning based curve fitting methods with square loss function, 
see    subsections \ref{ssect:efm} and \ref{ssect:log_normal}.


We then proceed with specific applications to Covid-19 scenarios.
Here we present, now using in addition data-driven machine learning based
curve fitting, some of our model's predictions to selected countries and US federal states, which are based on the currently existing respective data for these locations provided by the most recent numbers 
supplied by the Johns Hopkins University Covid-19 database.

\

Before mentioning and commenting upon other related works,
let us adopt from now on, and throughout all following parts 
of the present work, the following \emph{convention:}
as we  shall consider only reported cases  in our paper,
 we will omit the adjective 'reported'  from  'reported    cases  of infected'.

 \color{black}
 

In \cite{Ber20}, the authors  present three basic 'macroscopic' models to fit  data emerging from local and national
governments:  exponential growth, self-exciting branching process and compartmental models. The compartmental  models  are the classical SIR and SEIR  models,   the  self-exciting branching process  has been used 
before with regard to treating Ebola disease outbreaks   and  other dynamics
of social interaction.  In the exponential  growth model,  
the  number $ I(t)$ of  infectious  individuals  at time $t$ 
is expressed  as  $I(t) =I_0  e ^{\alpha t}$,  where $ \alpha$  is the rate constant.  The   exponential growth  model is  related to our approach, in which the exponential function is modelling the reported infections.  However, since this is a one-parameter model, 
it works only well for fitting the data at the beginning of an outbreak.

In \cite{NI} the authors  use a log-normal density function with three parameters to fit the daily reported cases. But since
they tried to fit the data with only one function, the curve of reported cases may not  be well fitted, 
since there are usually several waves of the epidemic for a period while one function presents only one wave. As explained above, our wavelet approach does overcome this difficulty.

In \cite{Tuli}, the authors use the function $f(x)=k. \gamma . \beta . \alpha^\beta . x^{-1-\beta} .  \exp(-\gamma (\alpha/x)^\beta)$ to fit all data. This method, too, can fit the data only for one wave.

In \cite{De} the authors fit the data of daily reported cases with a two-wave model, using the sum of two Gaussian functions.

In \cite{Chowell}, the authors introduce an  epidemic model 
composed of
overlapping sub-epidemic waves, where each single one is  
a generalized logistic growth model given
by solution of differential
equations. A short term forecast of the Covid-19 epidemic in China from February 5th to February 24th was given in \cite{Roo} using  three phenomenological models (generalized logistic growth model, the Richards growth model,  sub-epidemic wave model in \cite{Chowell}) and ensemble methods (see also \cite{Chowell20} for the ensemble approach in forecasting epidemic trajectories).

In \cite{KN},  a multi-wave model combining  several SIR models, namely,  a Multiple-Wave Forced-SIR model, was introduced to  fit the data of daily cases.

 Recently, Krantz et al. \cite{KPR} have proposed an approach to construct epidemic growth  models 
 using  {\it fractional} wavelets. These are  built   from the number of reported cases   to   construct    wavelets    that model   the  dynamics of the number of completed  cases (\cite{KPR}).  In their paper, the number  of completed  cases   is  the sum of the number of reported  cases  and   the number  of unreported  cases. Furthermore, the proposed approach there is  to update their models assuming the  availability of    the    reporting  error   which  improves  over  time  and tends  to zero       eventually.  This  assumption appears to us, however, as a too   idealistic one.

Those two last  approaches are the ones which are most closely related to our own. However, while those use single waves  coming from solutions differential equations,  we use general wavelet functions such as  Gaussian functions, log-normal functions, Gompertz density functions, and Beta prime density functions,  which   all  satisfy our general condition of  being epidemic-fitted in the sense of Definition \ref{def:ef}.

 We also refer to \cite{GLW, Hao,HCMK,YZW, Her,Xue,RBC,ACG, Kuc,DGM,Saq,Man, ISNG,ASV, KXL, HV, WAV, XWY, LMSW1,LMSW2,LMW,CNM,SDR,SVD} for other approaches on modelling and  forecasting the spread of Covid-19 epidemic using deep learning, machine learning, time series analysis, network model,  stochastic model and deterministic compartmental framework.

\

\noindent
The remaining parts of the present paper are  organized as follows:  

\

In section  \ref{sect:epic_model}  we   first  recall the  notion of a wavelet (Definition \ref{def:wavelet}) and the
fundamental theorem  of wavelet theory (Theorem \ref{thm:dau}) whiche we are going to put to use in the sequel. We proceed by introducing  
the  notion of   an epidemic-fitted  (EF) wavelet (Definition \ref{def:ef})  and   propose  our  method for  modelling
epidemic dynamics  (Proposition \ref{prop:ansatz}), justified by  the   fundamental theorem \ref{thm:dau}.   
In  Section \ref{sec:log} we consider  several  important  examples  of EF wavelets and     impose    constraints on an EF wavelet  to be   suitable as a basic  EF wavelet in  epidemic  dynamics. 
In Section \ref{sect:numerical}  we present applications  of   our   method  to  modelling and forecasting  the current spread of
Covid-19  in France, Germany, Italy, the Czech Republic and several US federal states, all based on the most recent JHU data.

\medskip
\noindent {\bf Acknowledgement.} It is our pleasure to thank
Tat Dat Tran and Vit  Fojtik for useful suggestions
and comments on an earlier version of this paper.

\section{Epidemic Modelling via Wavelet Theory and Machine Learning} \label{sect:epic_model}

\subsection{Wavelets}
In this subsection, we recall and collect some basic concepts and  facts from Wavelet Theory (cf. \cite{Dau,Mey96,Mey97}),
which will be needed in our approach for modelling epidemic dynamics.
\begin{definition}\label{def:wavelet}\cite[p.24]{Dau}  
 A {\it wavelet}  or (mother wavelet) is a function $\psi\in L^1(\mathbb{R})$ such that the following admissibility condition holds:
 \begin{equation}\label{eq:cond}
    C_\psi=  \int_{-\infty}^\infty|\hat \psi(\xi)|^2\frac{d\xi}{|\xi|}<\infty,
 \end{equation}
 where $\hat \psi $ is the Fourier transform of $\psi$, i.e,
 $\hat\psi(\xi) =\int_{\mathbb{R}} \psi(x) e^{-i\xi x} dx .$
\end{definition}
 
 Notice that condition \eqref{eq:cond}
 is only  satisfied if $ \hat \psi(0)=0$ or $\int \psi(x) dx =0$. Conversely, we have the following sufficient condition for \eqref{eq:cond}.
 \begin{lemma}\cite[p.24]{Dau} \label{Lem:int_cond}
Let $\psi\in L^1(\mathbb{R})$ and $\int_\mathbb{R} \psi(x) dx =0$. If  $\int_\mathbb{R} |\psi(x)| (1+|x|)^\alpha dx<\infty$ for some $\alpha>0$,  then $|\hat \psi(\xi) |\leq C|\xi|^{\min(\alpha,1)}$ and $C_\psi<\infty$.  
 \end{lemma} 
  A basic example of  a wavelet is the function 
  $$\psi(t)=\frac{\sin(2\pi t)-sin(\pi t)}{\pi t}.$$
From a mother wavelet one can generate other wavelets (called {\it children wavelets}), using affine transformations (i.e., dilations and translations):
 
 $$ \psi_{a,b}(t)= \frac{1}{ \sqrt{|a|}} \psi \left(\frac{t-b}{a} \right), \quad (a,b)\in \mathbb{R}\times \mathbb{R}. $$

 These wavelets  provide  us with the following
 decomposition of $L^2(\mathbb{R})$:
 
 \begin{theorem}\label{thm:dau}\cite[Proposition 2.4.1 and  p. 25-26]{Dau}
 Let $\psi$ be a mother wavelet. Then any $f\in L^2(\mathbb{R})$ decomposes as 
 
 \begin{equation}
     f=C_\psi^{-1}\int_{\mathbb{R}^2} <f, \psi_{a,b}> \psi_{a,b}\frac{dadb}{a^2}, 
 \end{equation}
 strongly in $L^2(\mathbb{R})$, where $<,>$ denotes the standard product in $L^2(\mathbb{R})$, i.e.,
 
 \begin{equation}
     \lim_{ A_1,A_2, B\rightarrow \infty
     }\| f- C_\psi^{-1}\int_{1/A_1\leq|a|\leq A_2, |b|\leq B} <f, \psi_{a,b}> \psi_{a,b}\frac{dadb}{a^2}\|_{L^2}=0. 
 \end{equation}
 \end{theorem}
  
 Any function  $f\in L^2(\mathbb{R})$ can then be written as a superposition of $\psi_{a_k,b_\ell}$, i.e., 
 
 $$ f(x)  = \sum_{k,l}  \alpha_{k,\ell} \psi_{a_k,b_\ell}(x).$$
  We refer to \cite{Dau} for more details on the analysis of discrete wavelet decomposition, and, especially, 
  for precise formulas for the coefficients $\alpha_{k,\ell}$. 
  
  \

\emph{
Notice that from a machine learning point of view, finding the $\alpha_{k,\ell},a_k,b_\ell$  can be thought of as a curve fitting  problem,
 and this is how we will combine wavelet theory and machine learning techniques
 in our approach to modelling epidemic dynamics.
 }

\subsection{Epidemic-fitted wavelets and   modelling}\label{ssect:efm}
As we already explained in the introduction, the time development of an epidemic features local as well as global wave-type phenomena.

This leads us to the concept of 
epidemic-fitted wavelets. 
Informally speaking, such a wavelet
is given by a positive real function
$W:\mathbb{R} \to \mathbb{R}^{>0}$,
whose value $W(t)$ at a given time $t$ describes the number of new infected cases in a homogeneous population with respect to an epidemic that occurs in one wave only and thus will satisfy some sort of homogeneous compartmental model (without network structure).

Since we are interested in the  daily    infected cases, we can assume that $W(t)$ is strictly positive but tends to 0 when $t$ tends
to $\pm \infty$. 
Setting $w(t) = \ln W(t)$ so that $W(t) = e^{w(t)}$, 
the (multiplicative) growth rate of $W$ is its log-derivative:
$$ \frac{\dot W(t)}{W}  =  \dot w(t).$$

We wish $W(t)$ to 'start' at $t=a$, 
(reach its) 'peak' at $t=\chi$, and 'stop' 
at $t=b$ ($a < \chi < b$). 
This is to say that $w(a) = w(b) = 0$, $\dot w(\chi) = 0$,
$\dot w(t) > 0$ for $t < \chi$ and $\dot w(t) < 0$ for $t > \chi$. 

\begin{definition}\label{def:ef}
Given an interval $ (a,b)\subset \mathbb{R}$, $a\geq 0$,
 an {\it epidemic-fitted wavelet} is a positive real function $\psi\in L^1((a,b),  \mathbb{R}^+ )$ such that $\psi$ has start-peak-stop behavior, i.e, $\psi$
 satisfies $\lim_{x\rightarrow a^+} \psi(x)= \lim_{x\rightarrow b^-} \psi(x) =0$, and $\psi$ admits its maximum at some point in $(a,b)$. 
\end{definition}

We can interpret $\psi$ as a wavelet 
$\tilde\psi$ in the of  Definition \ref{def:wavelet} by simply setting
$\tilde \psi(x):=\psi(x)$ for $x\in (a,b)$,  $\tilde \psi(x):=-\psi(|x|)$
for $x< -a$, and $\psi(x)=0$ otherwise. Indeed, this definition implies that $\int_{\mathbb{R}} \tilde\psi(x)dx =0$ and  $ \int_{\mathbb{R}} | \tilde \psi(x)|(1+|x|)dx <\infty $,  hence $C_{\tilde \psi}<\infty$ by Lemma \ref{Lem:int_cond} and $\tilde \psi$ is a wavelet. 

First examples of EF wavelets
which come to mind are polynomial
functions of degree $3$
(restricted to some finite interval).   
Other examples of functions with start-peak-stop behavior are 
Gaussian functions, log-normal functions, Gompertz density functions
\begin{equation}
    \psi_{b,c}(x)= bc\exp(c+bx-ce^{bx}),
\end{equation}
and, in SIR models, 
the solution function giving
the number of $I(t)$, the number of infectious individuals. 
(cf. \cite{BST}), etc.

\medskip
In our applications to real data (see Section \ref{sect:numerical}),
we will
employ log-normal functions
as epidemic-fitted (EF) wavelets.

\medskip

For treating an epidemic, 
we will concentrate on the curve of daily  (reported)  infected cases, denoted by $RC(t)$, and try to understand the epidemic growth based on this information.   

 Theorem  \ref{thm:dau}  implies  that  
 our following  ansatz  is  'asymptotically' correct, 
 as the number $N$ grows to  infinity.
 In particular, numerical simulations
 involving bigger and bigger numbers $N$
 will lead to better and better accuracy.

 \begin{proposition}[Ansatz] \label{prop:ansatz} A positive function
 (or curve)  whose value is the number of infected cases at time $t$ is 
 representable as a finite linear combination of epidemic-fitted wavelets:

\begin{equation}\label{eq:series}
    RC(t) = \sum_{i=1}^N \alpha_i W_i(t,\theta_i),
\end{equation}
where
each such wavelet $W_i$ can be obtained from a basic (mother) EF wavelet $\psi$ by adding some parameters
$\theta_i=(\theta_{i}^1, \ldots , \theta_i^k)$.   

\end{proposition} 

Using   this ansatz,  we shall model epidemic  dynamics  by finding 
the wavelet series  coefficients  $\alpha_i$  and $\theta_i$  in  the  decomposition (\ref{eq:series}), 
when given the  number of  infected  cases over   a sufficient  long time  frame.   

This amounts to solving  a {\it curve fitting problem} in {\bf machine learning}.


\section{Epidemic-fitted (EF) wavelets}\label{sec:log}
In this section, we introduce some epidemic models with different basic (mother) 
epidemic-fitted (EF) wavelets. 
In section \ref{sect:numerical} we show  by fitting the Covid-19 data  that log-normal EF wavelet models are highly
compatible with the data and lead to very
good forecast projections.

\subsection{Gaussian EF wavelets}
 The standard Gaussian function is a   fundamental
 example of a function which has start-peak-stop behavior and exponential growth:

\begin{eqnarray*}
    \psi: \mathbb{R}^+ &\rightarrow & (0,1]\\
    x &\mapsto & \exp(-x^2/2).
\end{eqnarray*}
After dilating and translating, we
obtain a general  Gaussian function

$$\psi_{b,c}(x)=\exp\left(-\frac{(x-b)^2}{2c^2}\right).$$

We remark that in general $\lim_{x\rightarrow-\infty}\psi_{b,c}(x) = 0$, but for certain $b,c>0$ we have $\psi_{b,c}(0) \ll 1$. In this case, we can simply set  $\tilde \psi(x) = \max(\psi_{b,c}(x)  -  \psi_{b,c}(0), 0) $ as the corresponding Gaussian EF wavelet.

In \cite{De} the authors fitted the data of daily reported cases with a two-wave model using the sum of two Gaussian functions. 
However, since these are symmetric 
with respect to the the vertical line $x=b$, this model may be not compatible with the curve of daily cases. We will explain this point in further detail 
in the next section.

\subsection{Log-normal EF wavelets}\label{ssect:log_normal}

We define here the log-normal function, which is a Gaussian function
in which the variable $x$ is interchanged
by $\log x$:

\begin{eqnarray*}
    \psi_{b,c}: \mathbb{R}^+ &\rightarrow & (0,1]\\
    x &\mapsto & \exp(-\frac{(\log x-b)^2}{2c^2}).
\end{eqnarray*}
We then define the corresponding log-normal wavelet by extending $$\psi_{b,c}(x)=-\exp(-\frac{(\log (-x)-b)^2}{2c^2}), \quad \text{ for } x<0.$$ 
Thus we can rewrite 

$$\psi_{b,c}(x)=sgn(x)\exp\left(-\frac{(\frac{1}{2}\log (x^2)-b)^2}{2c^2}\right).$$

By dilating and translating, we
obtain a general log-normal EF  wavelet 

$$\psi_{b,c,d}(x)=\exp\left(-\frac{(\frac{1}{2}\log (x-d)^2-b)^2}{2c^2}\right).$$

\begin{figure}[!htb]
\centering
\includegraphics[width=0.7\textwidth]{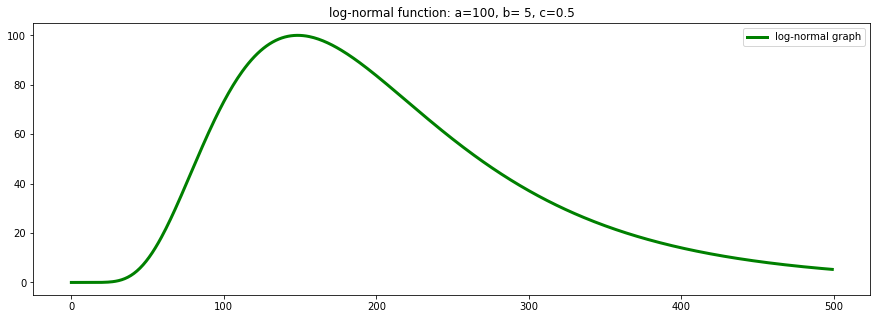}

\caption{Log-normal graph} \label{fig:log_nm}
\end{figure}

Figure \ref{fig:log_nm} depicts 
the graph of the log-normal function with scaling coefficient
$$\psi(x)=a\exp(-\frac{(\log x-b)^2}{2c^2}), x>0.$$

\subsection{Further examples of EF wavelets}
Based on probability distributions,
we can also choose many other functions 
to build a basic EF wavelet. 

For example, one can start here 
from Gompertz density functions
\begin{equation}
    \psi_{b,c}(x)= bc\exp(c+bx-ce^{bx}),
\end{equation}
or
Beta prime density functions
\begin{equation}
    \psi_{b,c}(x)= x^{b-1}(1+x)^{-b-c}/B(b,c),
\end{equation}
where $B$ is the Beta function. 

For appropiately chosen parameters $b,c$, they all satisfy the epidemic-fitted condition in Definition \ref{def:ef}

Another important class of EF wavelets is given by the function reporting
the  number of infectious individuals $I(t)$ in compartmental SIR models
and their variations (such as SEIR
and SIRD models, etc.) 

The SIR (compartmental) model  was introduced by  W. O. Kermack and A. G. McKendrick \cite{KM}, in which they considered a fixed population with only three compartments, and the numbers S(t)
(for 'susceptible'), I(t) (for 'infectious'),  and  R(t) (for 'recovered'
(or 'removed')).  

\begin{eqnarray}\label{eq:sir}
\frac{dS}{dt}&=&-\frac{\beta IS}{N}\\
\frac{dI}{dt}&=&\frac{\beta IS}{N}-\gamma I\\
\frac{dR}{dt}&=&\gamma I.
\end{eqnarray}

\begin{figure}[!htb]
\centering
\includegraphics[width=0.6\textwidth]{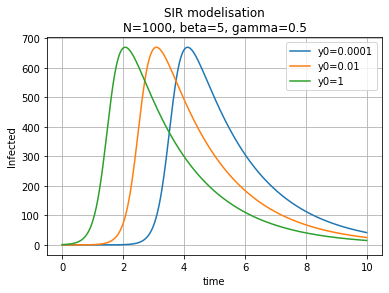}

\caption{Infectious individuals $I(t)$ for different initial conditions} \label{fig:sir}
\end{figure}

In Figure \ref{fig:sir} these curves show the number of infectious individuals  $I(t)$.

In general, $I(t)$ is an implicit function defined by a system of differential equations, which can lead to difficulties when trying to fit the data. However, 
we can here use the implicit solutions for simple SIR models which were deduced recently in \cite{BST}.

\subsection{Choosing suitable EF wavelets}\label{ssect:choose}
We explain here  how to choose good EF wavelets for building an epidemic model. 

The first  criterion to meet  is the start-peak-stop behavior as discussed in Section \ref{sect:epic_model}.

Our second criterion  is based on the following analysis of the number $I(t)$ of infectious individuals in SIR model:
\begin{eqnarray}
\frac{dS}{dt}&=&-\frac{\beta IS}{N}\\
\frac{dI}{dt}&=&\frac{\beta IS}{N}-\gamma I\\
\frac{dR}{dt}&=&\gamma I.
\end{eqnarray}

A closer look at SIR models reveals that  the    number $S(t)$ of susceptible individuals  is decreasing in time.
Therefore the number $I(t)$ of infectious
also grows less and the rate of infectious, i.e., $dI/dt$, before the peak is always less than the one after the peak. This is an important
criterion when choosing EF wavelets.  

{\it Log-normal EF wavelets actually turn out
to be very good candidates in this regard:}

Indeed, the first good point here is the  start-peak-stop behavior, where the start for a log-normal wavelet is at $x=0$ (or near 0),  the peak  is achieved at $x=e^b$ and the stop depends on the constant $c$. The second good one is that at the same value of $\psi$, the rate of  the curve before the peak is less than the one after the peak. 

This can be easily seen as follows. 
The derivative of $\psi_{b,c}$ is
\begin{equation}
    \psi'_{b,c}(x)=\psi_{b,c}(x)\frac{-(\log x-b)}{c^2 x}.
\end{equation}
Now suppose that $\psi(x_1)=\psi(x_2)$ with $x_1<e^b<x_2$, then $ |\log x_1-b|=|\log x_2-b|$. Therefore we have
$$ |\psi'_{b,c}(x_1)|= \frac{x_2}{x_1}|\psi'_{b,c}(x_2)|< |\psi'_{b,c}(x_2)|$$
as required.

These are the main reasons why we first chose log normal functions as basic EF wavelets for our numerical 
simulations (see   Section \ref{sect:numerical}).

We also remark that in \cite{NI} the authors used the log-normal density function, i.e.,  $f_{a,b,c}(x)=
\frac{a}{\sqrt{2\pi } c x}\psi_{b,c}(x)$, to fit   the number of daily reported cases. However, since they used only one single function, and since there are in general many waves of the epidemic,
the data  may not be  well-fitted enough
to produce realistic projections.

\section{Data-driven Numerical Forecasts} \label{sect:numerical}

In this section,  using log-normal 
EF wavelets we provide numerical results on the fitting  and forecasting  of daily new cases of Covid-19 epidemic  for some European countries and US federal states. 


\subsection{The log-normal wavelet model}\label{ssect:log_normal}
 Our EF wavelet model for the curve of daily new cases is a finite representation by log-normal EF  wavelest introduced in Section \ref{ssect:log_normal}:

$$W(t)= \sum_{i=1}^N  a_i \psi_{b_i,c_i}(t),$$
where $a_i,b_i,c_i$ are parameters, $N$ is the number of  log-normal EF  wavelets, and $t$ is the time variable. 

We intend to find the parameters $a_i,b_i,c_i$ such that $W(t)$ is close to the daily infected number $RC(t)$ by a  suitable loss function $L(\cdot,\cdot)$.  In other words, we want to find parameters which minimize $L(W,RC)$.

 For our numerical simulations 
 presented in the next section of this work,
 we shall use the Levenberg–Marquardt algorithm (cf. \cite{Lev, Mar}) for the least squares  loss function.
 The main advantage of this approach is
 that the loss function helps us to force the peaks of EF wavelets close to the peaks of real data. 
 
 \medskip
The number of log-normal wavelets $N$ depends on the data of each population level, since it presents the numbers sub-epidemic.  In our numerical simulations, we first try with $N=3, 5$. 
It would be interesting to estimate $N$ before fitting the model. 
Otherwise, we will need to choose $N$ sufficiently large, and redundant wavelets will have very small coefficients and,
correspondingly, very little effect.

\subsection{Data and smoothing} 
We will be using n the data supplied
by the 
Johns Hopkins University Center \cite{JH},
noting, however, that
almost all data from countries or US federal states  are subject to (high)  noise. One of the main reason for this is the reporting delay (cf. \cite{CS,Har}).

As explained in \cite{Har}, {\it ``there will be two main sources of
delay in monitoring trends. First of all, there will be a testing delay between the actual date when an
individual becomes infected and the date when that individual is ultimately tested. Second,
unless test samples are very rapidly processed, there will be a further reporting delay between
the date of testing and the date the test results are communicated by the reporting entity"}.

In order to reduce noise,  
 we  do smooth out  the real data using a (two-sided) moving average method (cf. \cite[Chapter 3]{MWH}, \cite{Hyn}, \cite{Sim}).
 A moving average is a time series constructed by taking averages of several sequential values
of another time series which is a type of mathematical convolution.  In statistics,  two-sided  moving averages are used to {\it smooth} a time series in order to estimate or highlight the underlying trend.  If we represent the original
time series by $x_1,...,x_n$, then a (simple) two-sided moving average of the time series will be  given by

$$\bar x_i=\frac{1}{2d+1} \sum_{k=i-d}^{i+d} x_k.$$

 If the data are showing a periodic fluctuation, moving averages of equal length period will eliminate the periodic variations (cf.  \cite{MWH,Hyn}). 
 Observing various population levels 
 indicates that there is  periodic fluctuation of 7 days on the data,  
 and hence  we will  take the average of 7 days
 $$\overline{RC}(i) = \frac{1}{7}\sum_{k=i-3}^{i+3} RC(k).$$


\subsection{Projections and validations for the Czech Republic, France, Germany, and Italy} 
\subsubsection{Projections from 25/10/2020}

\begin{figure}[htbp]
\centering
\includegraphics[width=0.8\textwidth]{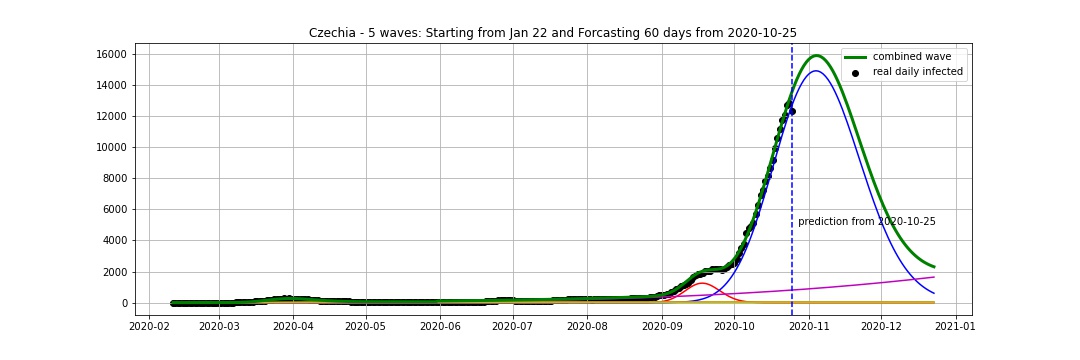}
\caption{Czechia: fitting and forecasting from 25/10 with 5 wavelets} \label{fig:Cz25}
\end{figure}

\begin{figure}[htbp]
\centering
\includegraphics[width=0.8\textwidth]{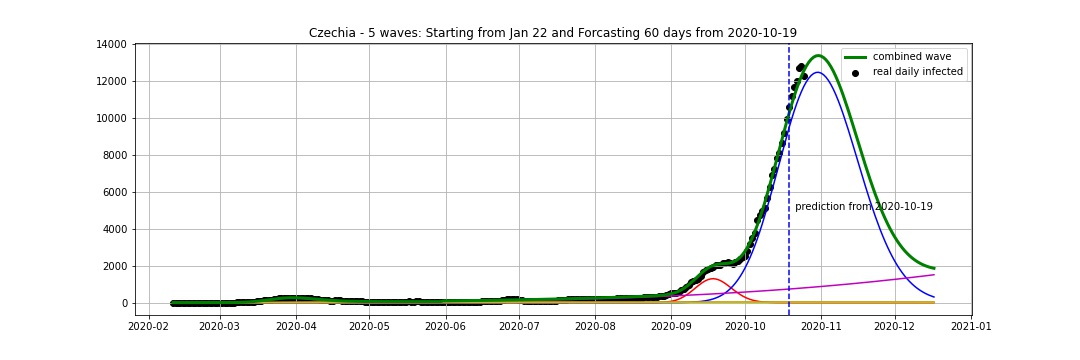}
\caption{Czechia: fitting and forecasting from 19/10 with 5 wavelets} \label{fig:Cz19}

\end{figure}

\begin{figure}[htbp]
\centering
\includegraphics[width=0.8\textwidth]{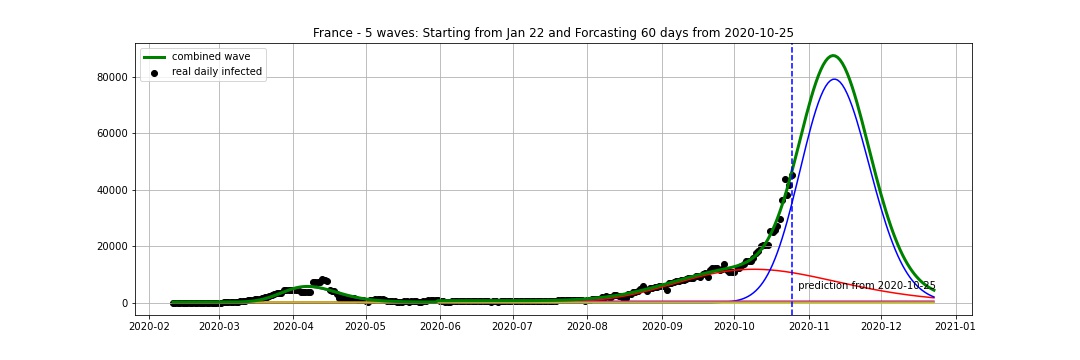}

\caption{France: fitting and forecasting from 25/10 with 5 wavelets.  Our model predicts a new wave starting from October 2020.} \label{fig:Fr25}
\end{figure}

\begin{figure}[htbp]
\centering
\includegraphics[width=0.8\textwidth]{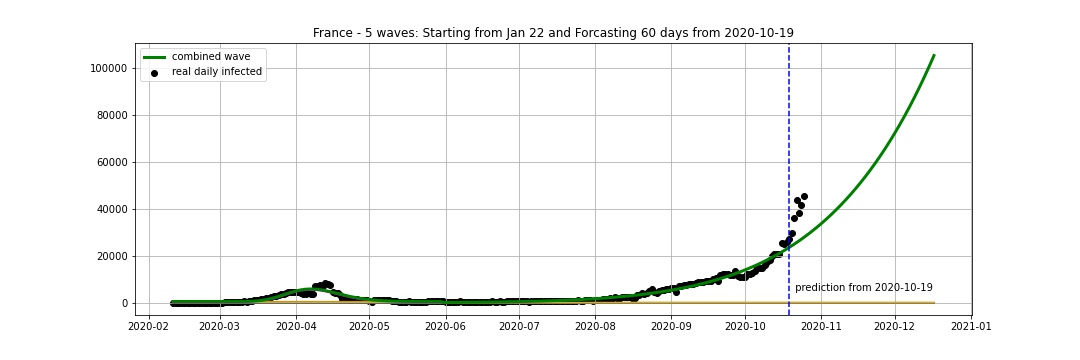}
\caption{France: fitting and forecasting from 19/10 with 5 wavelets.} \label{fig:Fr19}

\end{figure}

\begin{figure}[htbp]
\centering
\includegraphics[width=0.8\textwidth]{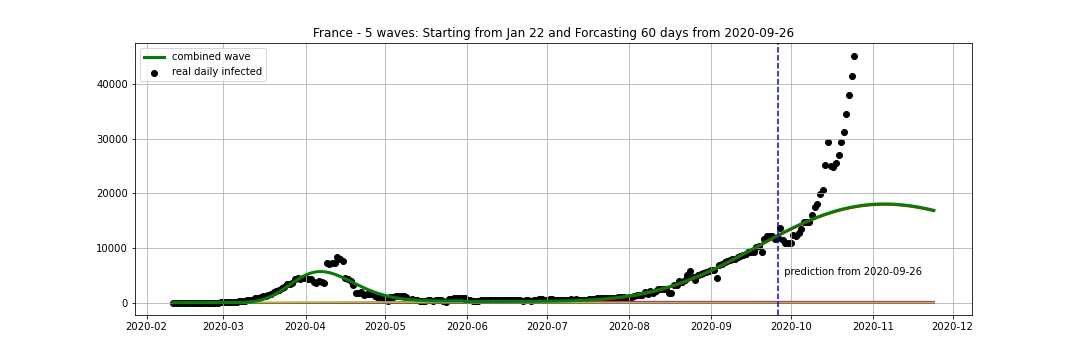}
\caption{France: fitting and forecasting from 26/09 with 5 wavelets} \label{fig:Fr26}

\end{figure}

\begin{figure}[htbp]
\centering
\includegraphics[width=0.8\textwidth]{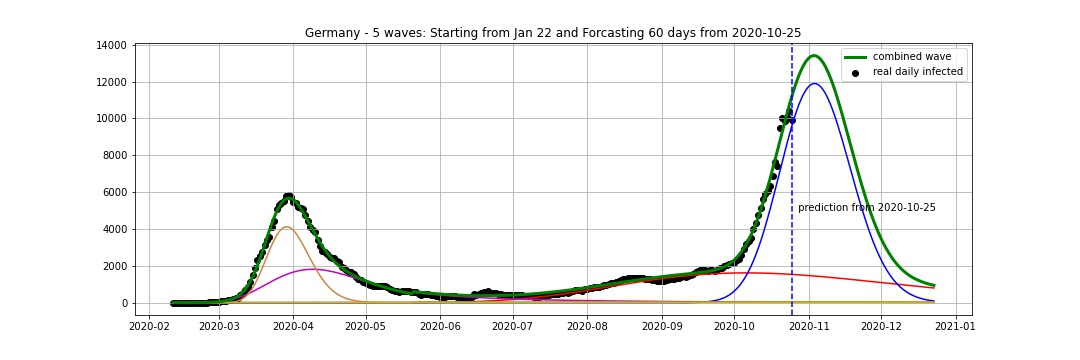}

\caption{Germany: fitting and forecasting from 25/10 with 5 wavelets} \label{fig:Ger25}
\end{figure}

\begin{figure}[htbp]
\centering
\includegraphics[width=0.8\textwidth]{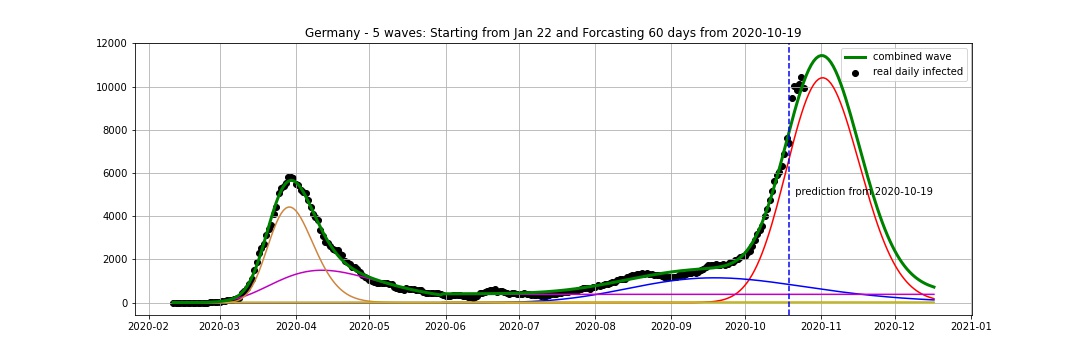}

\caption{Germany: fitting and forecasting from 19/10 with 5 wavelets} \label{fig:Ger19}
\end{figure}

\begin{figure}[htbp]
\centering
\includegraphics[width=0.8\textwidth]{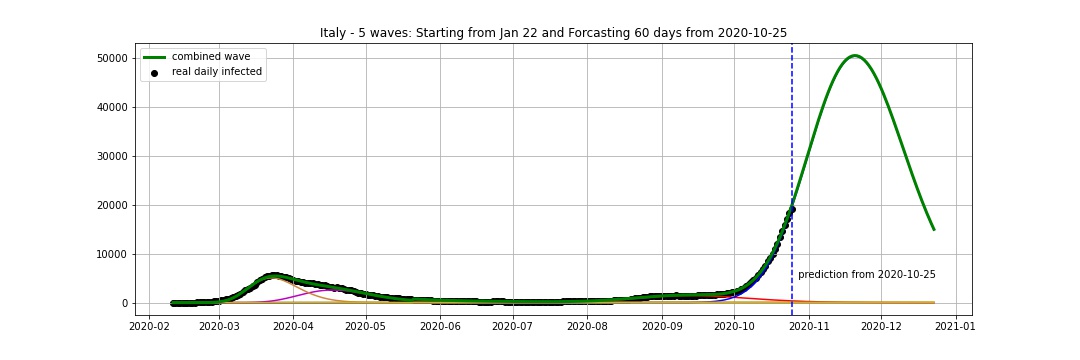}

\caption{Italy: fitting and forecasting from 25/10 with 5 wavelets} \label{fig:It25}
\end{figure}

\begin{figure}[htbp]
\centering
\includegraphics[width=0.8\textwidth]{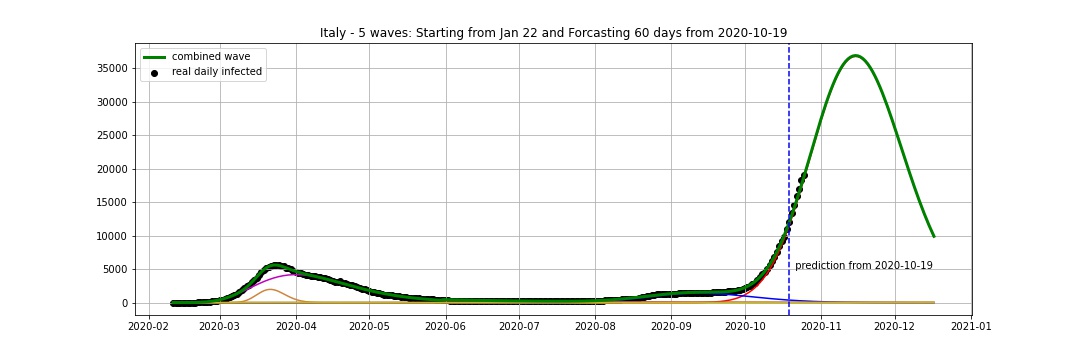}

\caption{Italy: fitting and forecasting from 19/10 with 5 wavelets} \label{fig:It19}
\end{figure}

In Figures \ref{fig:Cz25}, \ref{fig:Cz19},  \ref{fig:Fr25}, \ref{fig:Fr19}, \ref{fig:Fr26},
\ref{fig:Ger25}, \ref{fig:Ger19},  \ref{fig:It25}, \ref{fig:It19},   the  green curve shows 
the approximate number of daily confirmed new cases and also a possible scenario with a 60-day projection for the Czech Republic (or, in short: Czechia), France, Germany, and Italy. 
Other curves present log-normal EF  wavelets where each one can be seen as a sub-epidemic, localised both in time and location. 
These EF  wavelets then give us  the nowcasting for the epidemic situation for each population level, i.e.,  forecasts  present sub-epidemics,  recent sub-epidemics and the combination of sub-epidemics.

\medskip
For validation, we use the metric {\it relative percentage difference:}
\begin{equation}
    {\rm err}_i=\frac{| y_i- \hat y_i |}{y_i},
\end{equation}
where $y_i$ is the real data  at day $i$ smoothed by $7$-days moving average, $\hat y_i$ is the prediction of our model.  

 We fit our model with the data of daily cases until 19/10 and  keep last 6 days (20-25/10) for the validation set, then get the average error of 4.17\% for Czechia, 7.48\% for Germany and  3.25\% for Italy.

\begin{center}
 \begin{tabular}{||c c c c c||} 
  \hline\hline
& & Czechia& & \\ [0.5ex] 
 \hline
 day & real data & smoothing &prediction & error \\ [0.5ex] 
 \hline\hline
 20/10 & 11984 & 11173 & 10730 & 3.96\% \\ 
 \hline
 21/10 & 14969 & 11710 &  11161 & 4.68\% \\
 \hline
 22/10 & 14150
  & 12030 & 11564 & 3.87\% \\
 \hline
 23/10 & 15258 & 12689 & 11934 & 5.95\% \\
 \hline
  24/10 & 12474  & 12830 &12269 & 4.37\% \\ 
 \hline
 25/10 & 7300  & 12295 & 12564 & 2.18\% \\ [1ex] 
 \hline
\end{tabular}

\medskip
\begin{tabular}{||c c c c c||} 
  \hline\hline
& &Germany & & \\ [0.5ex] 
 \hline
 day & real data & smoothing& prediction & error \\ [0.5ex] 
 \hline\hline
 20/10 & 8523 & 9472 &  8346 &  11.88\% \\ 
 \hline
 21/10 & 12331 & 10019 & 8763& 12.53\% \\
 \hline
 22/10 & 5952
  & 9861 & 9164 & 7.06\% \\
 \hline
 23/10 & 22236 & 10105 &9545 & 5.54\% \\
 \hline
  24/10 & 8688  & 10421 &9902 & 4.98\% \\ [1ex] 
  \hline
   25/10 & 2900  & 9944 & 10231& 2.88\% \\ [1ex] 
 \hline
\end{tabular}

\medskip
\begin{tabular}{||c c c c c||} 
  \hline\hline
& &Italy & & \\ [0.5ex] 
 \hline
 day & real data & smoothing& prediction & error \\ [0.5ex] 
 \hline\hline
 20/10 & 10871 & 13322 &  13000 & 2.41\% \\ 
 \hline
 21/10 & 15199 & 14567 & 14080& 3.34\% \\
 \hline
 22/10 & 16078
  & 15934 & 15203 &  4.58\% \\
 \hline
 23/10 & 19143 & 17034 &16364 & 3.93\% \\
 \hline
  24/10 & 19640  & 18266 &17557 & 3.88\% \\ [1ex] 
  \hline
   25/10 & 21273  & 19033 & 18777& 1.34\% \\ [1ex] 
 \hline
\end{tabular}

\end{center}

 However, we get an average error of 32.61\% for France (see Figure \ref{fig:Fr19}) on the validation set from 20-25/10. We remark here that in some  3 consecutive days the total cases of France is constant in data  Johns Hopkins University \cite{JH}, and the total cases are updated by summing up for the day after these 3 days. For example 09-11/10 have same  732434 total cases, 16-18/10  have same  876342 total cases. This makes the daily reported cases are zero in some 2 consecutive days. Using moving average of 7 days we overcome this situation and then use the smoothing data for the projections shown in Figures \ref{fig:Fr25}, \ref{fig:Fr19}, \ref{fig:Fr26}.

\newpage

\subsubsection{Updated projections from 09/11/2020}
\begin{figure}[htbp]
\centering
\includegraphics[width=0.8\textwidth]{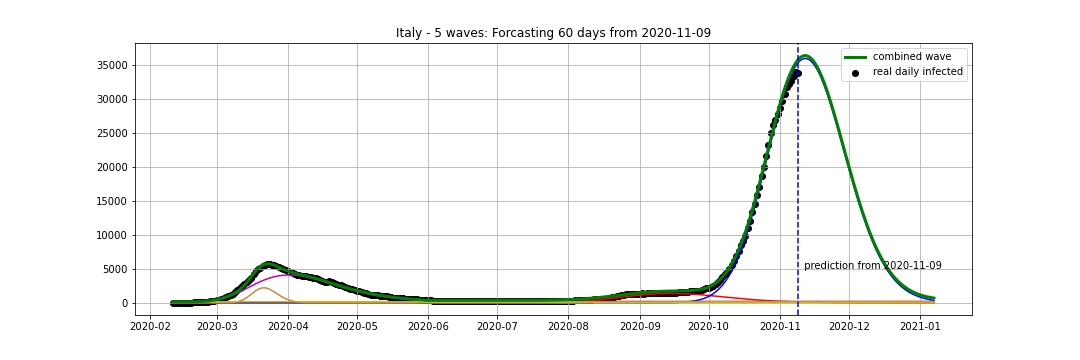}

\caption{Italy: fitting and forecasting from 09/11 with 5 wavelets} \label{fig:It_11}
\end{figure}

\begin{figure}[htbp]
\centering
\includegraphics[width=0.8\textwidth]{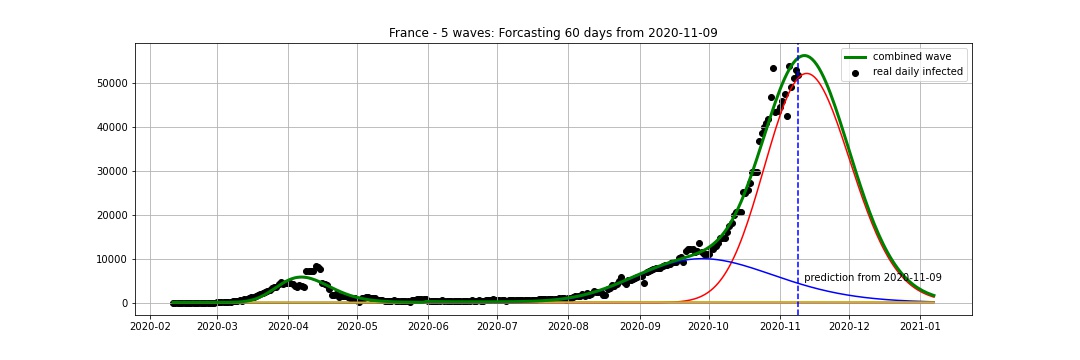}

\caption{France: fitting and forecasting from 09/11 with 5 wavelets} \label{fig:fr_11}
\end{figure}

\begin{figure}[htbp]
\centering
\includegraphics[width=0.8\textwidth]{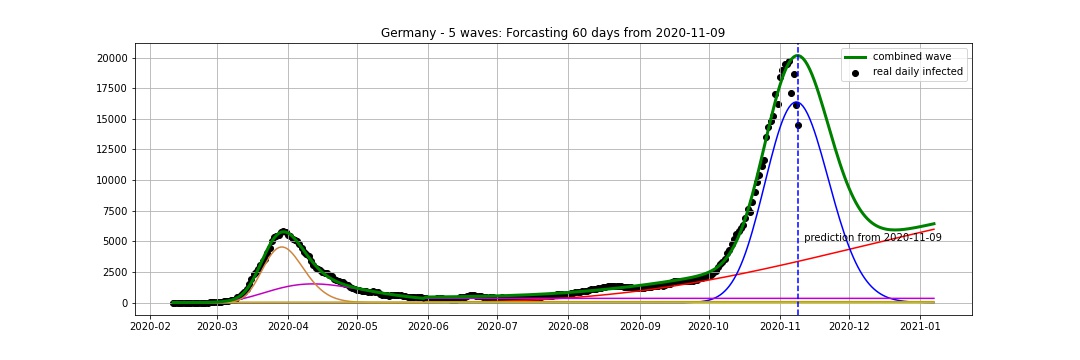}

\caption{Germany: fitting and forecasting from 09/11 with 5 wavelets} \label{fig:ger_11}
\end{figure}

\begin{figure}[htbp]
\centering
\includegraphics[width=0.8\textwidth]{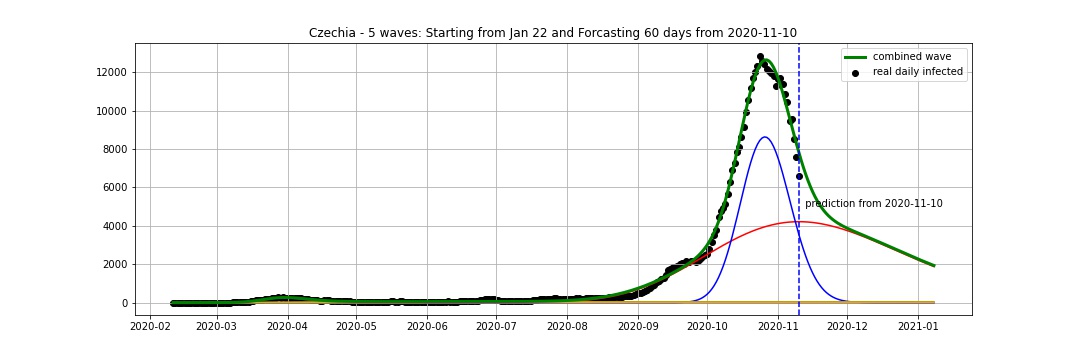}

\caption{Czechia: fitting and forecasting from 10/11 with 5 wavelets} \label{fig:cz_11}
\end{figure}

Figures \ref{fig:It_11}, \ref{fig:fr_11}, \ref{fig:ger_11}, \ref{fig:cz_11} show the projection from November 09, 2020. 

\newpage
\subsection{Projections for federal states in the United States}

\begin{figure}[htbp]
\centering
\includegraphics[width=0.8\textwidth]{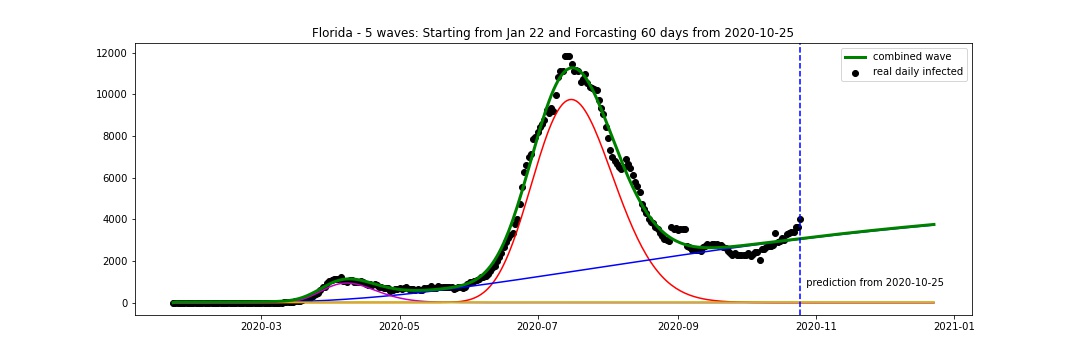}
\caption{Florida: fitting and forecasting from  25/10} \label{fig:flo}
\end{figure}

\begin{figure}[htbp]
\centering
\includegraphics[width=0.8\textwidth]{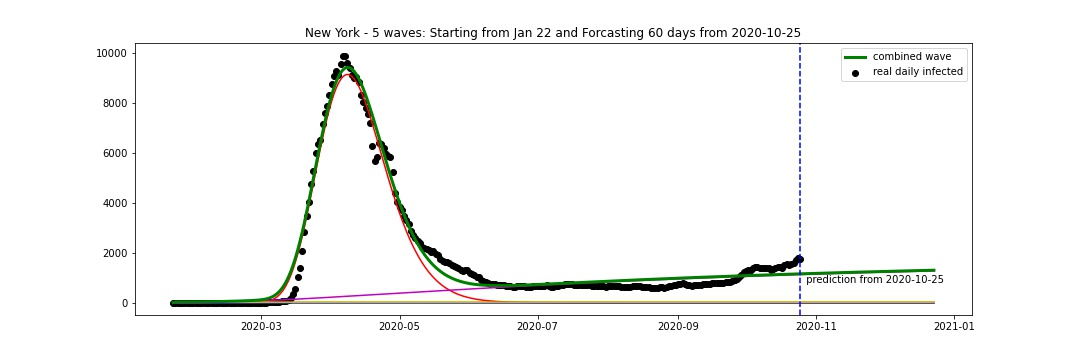}
\caption{New York: fitting and forecasting from  25/10} \label{fig:ny}

\end{figure}

In Figures \ref{fig:flo},   \ref{fig:ny}, 
   the  green curve shows the projection for Florida, New York, from  25/10/2020.   

\medskip
\noindent
{\bf Updated projections for Florida and New York from 10/11/2020}

\begin{figure}[htbp]
\centering
\includegraphics[width=0.8\textwidth]{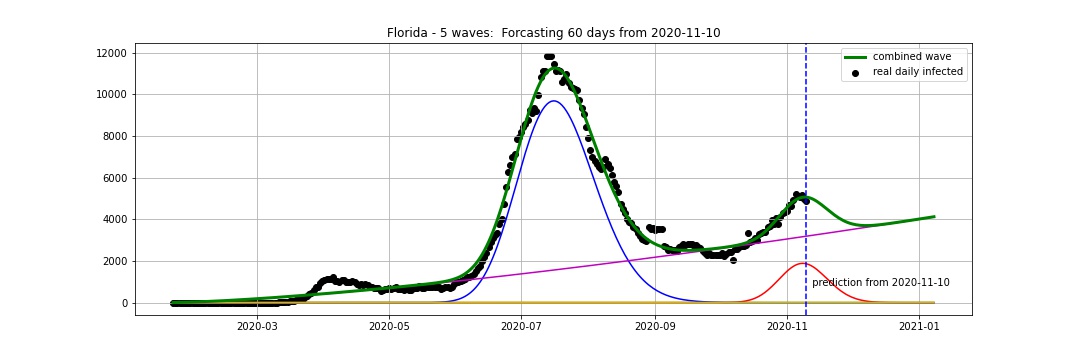}
\caption{Florida: fitting and forecasting from  10/11/2020} \label{fig:fl_2}
\end{figure}

\begin{figure}[htbp]
\centering
\includegraphics[width=0.8\textwidth]{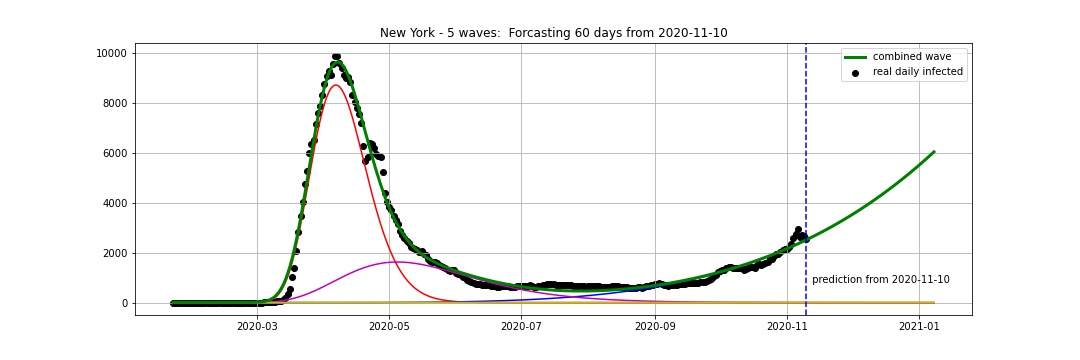}
\caption{New York: fitting and forecasting from  10/11/2020} \label{fig:ny_2}

\end{figure}

Figures \ref{fig:fl_2} and  \ref{fig:ny_2}  
show  the projection for Florida, New York, from  November 10, 2020.   

\newpage
\section{Conclusion and Outlook}\label{sec:concl}
\

(1) The numerical  results in the  last section of  our paper suggest that our  models are  actually able to predict  many  days  ahead  the  number  of  daily infected Covid-19 individuals in many different  countries. In particular, our approach also gives reasonable results for the epidemic situation on population levels by precising  sub-epidemics corresponding to EF wavelets.  

(2) For solving the curve-fitting problem in our  model  selection, we only have to use relatively few parameters. The model  can be seen as a neural network containing only one hidden layer  with a log-normal function  activation, entailing that we do not  have to deal with over-fitting problems and  that the estimation error of our model is low \cite[p.65]{SSBD2014}.

(3) Our method of modelling  of   the   number of daily reported  cases of  infectious individuals does also apply to  other  epidemics characteristics, e.g., to the number of active cases - and thus is also important  for health care system  decisions. 

(4) In future work, we will present
refinements of our approach as well as
refinements of the curve fitting techniques
employed here.


(5) In further future work, we will extend our approach based on epidemic-fitted wavelet approach
to situations where EF wavelets  are multivariate functions of time variables, measurement levels, or other variables such as death rate, recovery rate, etc.

{\small

\medskip
T\^o Tat Dat \\
Centre de Mathématiques Laurent-Schwartz, 
École polytechnique, 
Cour Vaneau, 91120 Palaiseau, France and  Institut de mathématiques de Jussieu – Paris Rive Gauche, Sorbonne Université, Campus Pierre et Marie Curie, 4, place Jussieu 
75252 Paris, France (current adress)

email: tat-dat.to@imj-prg.fr

Protin Frédéric,  Nguyen T.T. Hang, Martel Jules, Nguyen Duc Thang,  Charles Piffault,  Figueroa Susely \\
Torus Actions SAS, 
3 Avenue Didier Daurat,  31400 Toulouse

email: [protin;hangntt;jules;ndthang;charles.piffault; fsusely]@torus-actions.fr

 Rodr\'iguez Willy \\
 Ecole Nationale de l'Aviation Civile, 7 Avenue Edouard Belin, 31400 Toulouse

email: willy.rodriguez@enac.fr

H\^ong  V\^an  L\^e\\
Institute  of Mathematics of the Czech Academy of Sciences,
Zitna 25, 11567  Praha 1, Czech Republic,
email: hvle@math.cas.cz

Wilderich Tuschmann\\
Fakult\"at f\"ur Mathematik,
Karlsruher Institut f\"ur Technologie (KIT),
Englerstr. 2,
D-76131 Karlsruhe, Germany, 
email: tuschmann@kit.edu

Nguyen Tien Zung\\
Institut de Mathematiques de Toulouse, Universit\'e Toulouse 3,
email: tienzung@math.univ-toulouse.fr
}


\begin{thebibliography}{99}
\bibitem{ASV}
Acuna-Zegarra M.A., 
Santana-CibriancdM. and Velasco-Hernandez X. J.
Modeling behavioral change and COVID-19 containment in Mexico: A trade-off between lockdown and compliance,  {\it Mathematical Biosciences.} {\bf 2020}. 325: 108370. \url{https://doi.org/10.1016/j.mbs.2020.108370}

\bibitem{ACG}
 Arenas A.,  Cota W., Gomez-Gardenes J., 
 Gomez S., 
 Granell C.,  Matamalas J.,  Soriano D.,  Steinegger B.
A mathematical model for the spatiotemporal epidemic spreading of COVID19.  {\bf 2020}, {\it MedRxiv}, doi 10.1101/2020.03.21.20040022. 

\bibitem{Bar1956}
Bartlett M. S.  Deterministic  and  stochastic models for recurrent  epidemics. {\it Berkeley Symp. on Math. Statist. and Prob.} Proc. Third Berkeley Symp. on Math. Statist. and Prob. (Univ. of Calif. Press, {\bf 1956}), 4, 81-109.
 \url{https://projecteuclid.org/euclid.bsmsp/1200502549}


\bibitem{Bar1957}
Bartlett M. S.   Measles periodicity  and community size. {\it   J. R. Stat. Soc. A}. {\bf1957}, 120, 48-70.

\bibitem{Ber20}
 Bertozzi A. L., et al. The challenges of modeling and forecasting the spread of COVID-19. {\it Proceedings of the National Academy of Sciences}. {\bf 2020}, 117 (29), 16732-16738. \url{https://doi.org/10.1073/pnas.2006520117}.
 
\bibitem{BDW2008}
Brauer F., van den Driessche P., Wu  J. (Eds.). {\it Mathematical epidemiology}.
Lecture Notes in Mathematics 1945,  Mathematical Biosciences Subseries.
Springer-Verlag, Berlin, 2008.


\bibitem{CS}
 Cavataio J. and Schnell  S.
Interpreting SARS-CoV-2 fatality rate estimates — A case for
introducing standardized reporting to improve communication. {\it  SSRN}. {\bf 2020} \url{https://dx.doi.org/10.2139/ssrn.3695733} 
\bibitem{BST}
Bohner M.,  Streipert  S.,  Torres D. F. M.
Exact solution to a dynamic SIR model. {\it Nonlinear Analysis: Hybrid Systems}. {\bf 2019}, 32,  228-238.


\bibitem{Chowell}
Chowell G, Tariq A, Hyman J. A novel sub-epidemic modeling framework for short-term forecasting epidemic waves. {\it BMC Med.} {\bf 2019}, 17, 164. \url{https://doi.org/10.1186/s12916-019-1406-6}.

\bibitem{Chowell20}
Chowell G., Luo R., Sun K., Roosa K., Tariq A., Viboud C.
Real-time forecasting of epidemic trajectories using computational dynamic ensembles, {\it Epidemics}, 
 {\bf 2020}, 30., 100379
\bibitem{CNM}
Cotta R.M., Naveira-Cotta C.P., Magal P. Mathematical Parameters of the COVID-19 Epidemic in Brazil and Evaluation of the Impact of Different Public Health Measures. {\it Biology} {\bf 2020}, 9, 220.
\bibitem{De}
De Noni Jr. A., et al. 
A two-wave epidemiological model
of COVID-19 outbreaks using MS-Excel, {\it medRxiv}. {\bf 2020},  doi: \url{https://doi.org/10.1101/2020.05.08.20095133}. 


\bibitem{Dau}
 Daubechies I. {\it Ten Lectures on Wavelets}. Society for Industrial and Applied Mathematics, 1992. 


\bibitem{DGM}
 Demongeot J.,  Griette Q. and  Magal P. SI epidemic model applied to COVID-19 data in mainland China, medRxiv, doi:  10.1101/2020.10.19.20214528

\bibitem{GLW} 
Guo Q.,  Li M., Wang C., et al. Host and infectivity 
prediction of Wuhan 2019 novel coronavirus using deep
learning algorithm. {\it bioRxiv}. {\bf 2020}, doi: \url{https://doi.org/10.1101/2020.01.21.914044}. 

\bibitem{Hao}
Hao Y., Xu T., Hu H., Wang P., Bai Y. {\it  
Prediction and analysis of Corona Virus Disease
2019}, {\it PLoS ONE}. {\bf 2020}, 15(10): e0239960. \url{https://doi.org/10.1371/journal.pone.0239960}. 

\bibitem{Har}
Harris J. E. {\it
Overcoming Reporting Delays Is Critical to Timely Epidemic Monitoring: The Case of COVID-19 in New York City},
{\it medRxiv}. {\bf 2020}.
doi: \url{https://doi.org/10.1101/2020.08.02.20159418}.


\bibitem{Her}
Hernand-Matamoros A., et al. 
Forecasting of COVID19 per regions using ARIMA models and polynomial functions. {\it Applied Soft Computing}. {\bf 2020}, 96, 106610.
\bibitem{HV}
Hernandez-Vargas E. A., Velasco-Hernandez, J. X.
In-host Mathematical Modelling of COVID-19 in Humans.  {\it Annual Reviews in Control} {\bf 2020}, to appear. doi:  10.1016/j.arcontrol.2020.09.006 


\bibitem{HCMK} Huang C.-Y., Chen Y.-H., Ma  Y.,  and Kuo  P.-H.
Multiple-Input Deep Convolutional Neural Network
2 Model for COVID-19 Forecasting in China, {\it medRxiv}. {\bf 2020}. doi:\url{https://doi.org/10.1101/2020.03.23.20041608}.



\bibitem{Hyn}
Hyndman R.J.   Moving Averages. In:  Lovric M. (eds.) {\it International Encyclopedia of Statistical Science.} Springer, Berlin, Heidelberg, 2011.  


\bibitem{ISNG}
 Iboi E., Sharomi O.,  Ngonghala C. and  Gumel A. B.
Mathematical Modeling and Analysis of COVID-19 pandemic in
Nigeria.  {\bf medRxiv}. {\bf 2020} . doi: 10.1101/2020.05.22.20110387


\bibitem{JH}
 Johns Hopkins University Center, Covid-19 data: \\  \hyperlink{https://github.com/CSSEGISandData/COVID-19}{https://github.com/CSSEGISandData/COVID-19}. 

\bibitem{KM}
Kermack W. O., McKendrick A. G.  A Contribution to the Mathematical Theory of Epidemics. {\it Proceedings of the Royal Society A.} {\bf 1927},115 (772), 700-721. 

\bibitem{KN}
 Kaxiras E.,  Neofotistos G.
Multiple Epidemic Wave Model of the COVID-19 Pandemic: Modeling Study.
{\it J Med Internet Res.} {\bf 2020},  22(7): e20912. doi:10.2196/20912.


\bibitem{KXL}
Kapoor A., Ben X., Liu L., Perozzi B., Barnes M., Blais M., O'Banion S.
Examining COVID-19 Forecasting using Spatio-Temporal Graph Neural Networks, {\it arXiv}. {\bf 2020}.  \url{https://arxiv.org/abs/2007.03113}

\bibitem{KPR}
Krantz P. P., Polyakov  P.,   Rao A. S .R. S.  
True epidemic growth construction through harmonic analysis. {\it Journal of Theoretical Biology}, {\bf 2020}, 494, 110243.

\bibitem{KR}
Keeling M. J.,   Rohani P.  {\it Modeling Infectious Diseases in Humans and Animals}. Princeton University  Press, 2008.

\bibitem{Kuc}
 Kucharski A. J. et al. Early dynamics of transmission and control of COVID-19: A mathematical modelling study. {\it Lancet Infect. Dis.} {\bf 2020} 20, 553–558 .

\bibitem{Lev}
Levenberg K. A Method for the Solution of Certain Non-Linear Problems in Least Squares. {\it Quarterly of Applied Mathematics.} {\bf 1944} 2 (2): 164–168.
 
\bibitem{LMSW1}
 Liu Z., Magal P., Seydi  O. and  Webb G. Predicting the cumulative number of cases for the COVID-19 epidemic in
China from early data. {\it Mathematical Biosciences and Engineering}. {\bf 2020} 17(4), 3040-3051.
\bibitem{LMSW2}
 Liu Z., Magal P., Seydi  O. and  Webb G.
 A COVID-19 epidemic model with latency period, {\it Infectious Disease Modelling} {\bf 2020} 5, Pages 323-337.
\bibitem{LMW}
 Liu Z., Magal P. and  Webb G. Predicting the number of reported and unreported cases for the COVID-19 epidemics in China, South Korea, Italy, France, Germany and United Kingdom. {\it Journal of Theoretical Biology}, {\bf 2021} 509, 21.
\bibitem{NI}
 Nishimoto1 Y.,  Inoue K.
Curve-fitting approach for COVID-19 data and its physical background, {\it medRxiv}. {\bf 2020}. doi: \url{https://doi.org/10.1101/2020.07.02.20144899}.

\bibitem{Man}
Manevski D. et al.
Modeling COVID-19 pandemic using Bayesian analysis with application to Slovene data. {\it Mathematical Biosciences.} {\bf 2020} 329.  108466 

\bibitem{MWH}
Makridakis S., Wheelwright S. C., Hyndman R. J.  {\it Forecasting: methods and applications, 3rd edition}. Wiley, New York, 1998.

\bibitem{Mar}
 Marquardt D. An Algorithm for Least-Squares Estimation of Nonlinear Parameters. {\it SIAM Journal on Applied Mathematics}. {\bf 1963} 11 (2): 431–441.

\bibitem{Mey96}
 Meyer Y.  and  Ryan D. {\it 
Wavelets: Algorithms and Applications}, Society for Industrial and Applied Mathematics, 1996.

\bibitem{Mey97}
 Meyer Y. {\it 
Wavelets, Vibrations and Scalings},  (CRM Monograph Series),
American Mathematical Society, 1997.
\bibitem{RBC}
Reiner, R.C., Barber, R.M., Collins, J.K. et al. Modeling COVID-19 scenarios for the United States. {\bf 2020} { \it  Nat Med }. doi: 10.1038/s41591-020-1132-9

\bibitem{Roo}
Roosa K.,  Lee Y.,  Luo R., Kirpich  A.,  Rothenberg R.,  Hyman J.M., Yan P.,  Chowell G. Real-time forecasts of the COVID-19 epidemic in China from February 5th to February 24th, 2020, {\it Infectious Disease Modelling}. {\bf 2020}, 5.  256-263
\bibitem{Saq}
Saqib, M. Forecasting COVID-19 outbreak progression using hybrid polynomial-Bayesian ridge regression model. {\it Appl Intell} {\bf 2020}. \url{https://doi.org/10.1007/s10489-020-01942-7}

\bibitem{SDR}
 Soubeyrand S.,  Demongeot J. and Roques L.
Towards unified and real-time analyses of outbreaks at country-level during
pandemics. {\it One Health}. {\bf 2020}, 100187.
\bibitem{SVD}
Seligmann H., Vuillerme N. and Demongeot J.
Summer COVID-19 third wave: faster high altitude spread suggests high UV adaptation.  {\it medRxiv} {\bf 2020}. DOI: 10.1101/2020.08.17.20176628
\bibitem{Sim}
Simonoff J. S.  {\it Smoothing Methods in Statistics, 2nd edition}. Springer-Verlag, New York, 1996. 

\bibitem{Sop}
Soper H. E. The interpretation of periodicity in disease prevalence. {\it J. Roy. Stat. Soc, Ser. A.} {\bf 1929}, 92, 34-61.

\bibitem{SSBD2014}
 Shalev-Shwartz S.,  Ben-David S. {\it Understanding Machine Learning: From Theory to Algorithms}.  Cambridge University Press, 2014.
 
\bibitem{Tuli}
Tuli S., et al. 
Predicting the growth and trend of COVID-19 pandemic using machine learning and cloud computing.{\it  Internet of Things}. {\bf 2020}, 11, 100222 \url{https://doi.org/10.1016/j.iot.2020.100222}.


\bibitem{Wang}   Wang J.  Mathematical models for COVID-19: applications, limitations, and potentials. {\it Journal of Public Health and Emergency}. {\bf 2020}, 4. doi: 10.21037/jphe-2020-05.  
\bibitem{Wan}
 Wang S. et al.
Modeling the viral dynamics of SARS-CoV-2 infection, {\it Mathematical Biosciences}. {\bf 2020} 328. 108438


\bibitem{WAV}
Wang L., Adiga A., Venkatramanan S., Chen J., Lewis B., Marathe M. Examining Deep Learning Models with Multiple Data Sources for COVID-19 Forecasting. {\it arXiv}. {\bf 2020, 10}. doi: \url{https://arxiv.org/abs/2010.14491}


\bibitem{Xue}
 Xue L. et al. 
A data-driven network model for the emerging COVID-19 epidemics in Wuhan, Toronto and Italy. {\it Math Biosci.}  {\bf 2020} 326.  108391

\bibitem{YZW} Yang Z.,  Zeng Z., Wang K., et al. Modified SEIR and 
AI prediction of the epidemics trend of COVID-19 in China under public health interventions. {\it J. of  Thoracic Disease}. {\bf 2020}, 12(3):165-174. doi: 10.21037/jtd.2020.02.64.

\bibitem{XWY}
Xiaoyong J., Yu-Xiang W. and Xifeng Y. Inter-Series Attention Model for COVID-19 Forecasting. {\it arXiv}. {\bf 2020, 10}. doi: \url{https://arxiv.org/abs/2010.13006}


\end{thebibliography}
\end{document}